\documentclass{amsart}
\usepackage[mathscr]{eucal}

\usepackage{xypic}
\usepackage{amsfonts}
\usepackage{amsmath}
\usepackage{amsthm}
\usepackage{amssymb}
\usepackage{latexsym}
\usepackage{eufrak}

\newtheorem{theorem}{Theorem}[section]
\newtheorem{lemma}[theorem]{Lemma}
\newtheorem{proposition}[theorem]{Proposition}
\newtheorem{corollary}[theorem]{Corollary} 
\theoremstyle{definition}

\newcommand{\Tr}{\text{Tr}}

\newcommand{\eps}{\varepsilon}

\newcommand{\h}{\mathfrak{h}}

\newcommand{\ben}{\begin{enumerate}}
\newcommand{\een}{\end{enumerate}}

\newcommand{\BZ}{{\mathbb Z}}

\newcommand{\BC}{{\mathbb C}}

\begin{document}

\title{On some representations 
of the rational Cherednik algebra}

\begin{abstract} 
We study lowest weight representations of the rational Cherednik
algebra attached to a complex reflection group $W$. 
In particular, we generalize a number of previous results 
due to Berest, Etingof, and Ginzburg. 
\end{abstract}

\author{Tatyana Chmutova}
\address{Department of Mathematics, Harvard University,
Cambridge, MA 02138, USA}
\email{chmutova@math.harvard.edu}

\author{Pavel Etingof}
\address{Department of Mathematics, Massachusetts Institute of Technology,
Cambridge, MA 02139, USA}
\email{etingof@math.mit.edu}

\maketitle

\section{Introduction}
In this paper we study the structure of some lowest weight
representations of the rational Cherednik algebra $H_c(W)$
attached to an irreducible complex reflection group $W$
(\cite{EG,GGOR}). 

The composition of the paper is as follows. 
In Section 2 we recall 
the definition and basic properties of the rational Cherednik algebra
and its representations. Then we describe the rank 1 case, and 
prove a few general results. In particular, we show that 
for real reflection groups, 
a finite dimensional quotient of the polynomial representation
is irreducible if and only if it is a Gorenstein algebra. 

In Section 3, we consider the special case $W=S_n$.   
Let $M_k$ be the polynomial representation of 
$H_k(S_n)$. Dunkl showed in \cite{Du} that if $r$ is a positive
integer not divisible by $n$ and $k=r/n$, then $M_k$ contains a copy 
of the reflection representation of $S_n$ in degree $r$ which consists of
singular vectors. Let $I_k$ be the $H_k(S_n)$-submodule in $M_k$ 
generated by these singular vectors. We compute the support of
the module $M_k/I_k$ as a $\mathbb C[\h]$-module. 
In particular, we show that the Gelfand-Kirillov dimension 
of $M_k/I_k$ is $d-1$, where $d$ is the greatest common divisor 
of $r$ and $n$. In the special case $d=1$, this 
implies that $M_k/I_k$ is finite dimensional. Using this fact and
the results of Section 2, we give a simple
proof of the result from \cite{BEG} that the module $M_k/I_k$ 
has dimension $r^{n-1}$ and is irreducible. 

In Section 4, consider the case when $W$ is 
the complex reflection group \linebreak
$S_n\ltimes (\Bbb Z/l\Bbb Z)^n$. We use a similar 
formula for singular vectors (due to Dunkl and Opdam, \cite{DO})
to study finite dimensional
representations of
$H_{c}(W)$. More specifically, for each 
positive integer $r$ not divisible by $l$, 
we define a hyperplane $E_r$ in the space of the functions $c$,
and for each $c\in E_r$ construct a quotient $Y_c$ of the polynomial
representation $M_c$ of dimension $r^n$, which is generically
irreducible. For $l=2$, we use the results of Section 2 to obtain
a more precise information about the set of $c\in E_r$ for which 
$Y_c$ is irreducible. 

{\bf Acknowledgments.} The research of P.E. was partially supported by
the NSF grant DMS-9988796, and was done in part for the Clay Mathematics 
Institute.

\section{Representations
of the rational Cherednik algebra for a general complex reflection group}

\subsection{Definitions and notation}
In this subsection we recall the standard theory 
of the rational Cherednik algebra \cite{GGOR,EG}.

Let $W$ be an irreducible complex reflection group with reflection 
representation ${\frak h}$ of dimension $\ell$. Let 
$c$ be a conjugation invariant complex function
on the set $S$ of complex reflections of $W$. 
The rational Cherednik algebra $H_c(W)$ 
is generated by 
$g\in W$, $x\in \h^*, y\in \h$, with defining relations
$$
gxg^{-1}=x^g, gyg^{-1}=y^g, [x,x']=0, [y,y']=0,
$$
$$
[y,x]=(y,x)-\sum_{s\in S} c_s (y,\alpha_s)(\alpha_s^\vee,x)s,
$$
for $x,x'\in \h^*$, $y,y'\in \h$, $g\in W$. Here $\alpha_s$ is 
a nonzero linear function on $\h$ 
vanishing on the reflection hyperplane for $s$, 
and $\alpha_s^\vee$ is a linear function on $\h^*$ with the same property, 
such that $(\alpha_s^\vee,\alpha_s)=2$. 

For any irreducible representation $\tau$ of $W$, 
let $M_c(\tau)$ be the standard representation 
of $H_c(W)$ with lowest weight 
$\tau$, i.e. $M_c(\tau)=H_c(W)\otimes_{\Bbb C[W]\ltimes \Bbb C[\h^*]}\tau$, 
where $\tau$ is the representation of 
$\Bbb CW\ltimes \Bbb C[\h^*]$, in which $y\in \h$
act by $0$. Thus, as a vector space, 
$M_c(\tau)$ is naturally identified with $\Bbb C[\h]\otimes \tau$. 

An important special case of $M_c(\tau)$ is 
$M_c=M_c(\Bbb C)$, the polynomial representation,
corresponding to the case when $\tau=\Bbb C$ is trivial. 
The polynomial representation can thus be 
naturally identified with $\Bbb C[\h]$. 
Elements of $W$ and $\h^*$ act in this space in the obvious way, 
while elements $y\in \h$ act by Dunkl operators
$$
\partial_y-\sum_{s\in S}\frac{2c_s}{1-\lambda_s}\frac{(\alpha_s,y)}{\alpha_s}
(1-s),
$$
where $\lambda_s$ is the nontrivial eigenvalue of $s$ in the 
dual reflection representation. 

An important element in $H_c(W)$ is the element 
$$
\bold h=\sum_i x_iy_i+\frac{\ell}{2}-\sum_{s\in S}\frac{2c_s}{1-\lambda_s}s,
$$
where $y_i$ is a basis of $\h$ and $x_i$ the dual basis of $\h^*$. 
This element is 
$W$-invariant and satisfies the equations $[\bold h,x]=x$ and $[\bold h,y]=-y$.
The category $\mathcal O$ of $H_c(W)$-modules is the category 
of $H_c(W)$-modules $V$, such that $V$ is the direct sum of finite 
dimensional generalized eigenspaces of $\bold h$, and the 
real part of the spectrum of $\bold h$ 
is bounded below. The standard representations $M_c(\tau)$ and their 
irreducible quotients $L_c(\tau)$ belong to ${\mathcal O}$. 
The character of 
a module $V\in {\mathcal O}$ is 
$\chi_V(g,t)=\Tr|_V(g t^\bold h)$, $g\in W$ (this is a series in $t$).  
For example, the character of $M_c(\tau)$ is 
$$
\chi_{M_c(\tau)}(g,t)=\frac{\chi_\tau(g)t^{h(\tau)}}{\det|_{\h^*}(1-gt)},
$$
where $h(\tau):=
\ell/2-\sum_s \frac{2c_s}
{1-\lambda_s}s|_\tau$ is the lowest eigenvalue of 
$\bold h$ in $M_c(\tau)$. 

We note that if $W$ is a real reflection group
then $\bold h$ can be included in an $sl_2$ triple
$\bold h,\bold E=\frac{1}{2}\sum x_i^2,\bold F=\frac{1}{2}\sum
y_i^2$, where $x_i$, $y_i$ are orthonormal bases of $\h^*$ and
$\h$, respectively (see e.g. \cite{BEG1}, Section 3).   

The module $L_c(\tau)$ can be characterized in terms of the 
contragredient standard modules. 
Namely, let $\hat M_c(\tau)=\tau^*\otimes_{\Bbb C[W]\ltimes \Bbb C[\h]}H_c(W)$
be a right $H_c(W)$-module, and $M_c(\tau)^\vee=\hat M_c(\tau)^*$ its 
restricted dual, which
may be called the contragredient standard module. 
Clearly, there is a natural morphism 
$\phi: M_c(\tau)\to M_c(\tau)^\vee$.
The module $L_c(\tau)$ is the image of $\phi$. 

Note that the map $\phi$ can be viewed as a 
bilinear form $B: \hat M_c(\tau)\otimes M_c(\tau)\to \Bbb C$. 
This form is analogous to the Shapovalov form in Lie theory. 

If $W$ is a real reflection group, then $H_c(W)$ admits a natural
antiinvolution \linebreak $x_i\to y_i$, $y_i\to x_i$, $g\to g^{-1}$ for
$g\in W$ (where $x_i,y_i$ are orthonormal bases of $\h^*$ and
$\h$ dual to each other). Under this antiinvolution, the right
module $\hat M_c(\tau)$ turns into the left module 
$M_c(\tau^*)$, so the form $B$ is a (possibly degenerate) pairing 
$M_c(\tau^*)\otimes M_c(\tau)\to \Bbb C$.  
Moreover, it is clear that if $Y,Y'$ are any quotients 
of $M_c(\tau),M_c(\tau^*)$ respectively, then $B$ descends to 
a pairing $Y'\otimes Y\to \Bbb C$
(nondegenerate iff $Y,Y'$ are irreducible). This pairing 
satisfies the contravariance equations $B(a,x_ib)=B(y_ia,b)$, 
$B(a,y_ib)=B(x_ia,b)$, and $B(ga,gb)=B(a,b)$. 

\subsection{The rank 1 case}
One of the main problems in representation theory 
of the rational Cherednik algebra is to compute the multiplicities 
of $L_c(\sigma)$ in $M_c(\tau)$ or, equivalently,
 the characters of $L_c(\sigma)$. This problem is hard and open in general. 
However, in the rank one case ($\ell=1$) this problem is 
trivial to solve. Nevertheless, it is an instructive example, and 
we will give the answer, omitting the proofs, which are straightforward.

In the rank $1$ case, $W=\Bbb Z/l\Bbb Z$, and 
the reflection representation is $\Bbb C$, with the generator $s$ of 
$W$ acting by $\varepsilon$, 
where $\varepsilon$ is a primitive root of unity of degree $l$. 
The function $c$ is a collection of numbers 
$(c_1,...,c_{l-1})$ (where $c_i=c_{s^i}$), and the 
algebra $H_c(W)$ is generated by three generators 
$x,y,s$ with defining relations 
$$
sx=\varepsilon^{-1}xs, sy=\varepsilon ys, s^l=1, 
$$
$$
[y,x]=1-2\sum c_js^j.
$$
The irreducible representations of $W$ are 
$\eta^j$, where $\eta(s)=\varepsilon$. 

Define the polynomial $f_c(z)=\sum \frac{2c_j}{1-\varepsilon^{-j}}z^j$.
The lowest eigenvalue of $\bold h$ on 
$M_c(\eta^j)$ is $\frac{1}{2}-f_c(\varepsilon^j)$.

\begin{theorem}\label{rank1} (i) 
The multiplicity of $L_c(\eta^m)$ in $M_c(\eta^p)$ 
is $1$ if $f_c(\varepsilon^p)-f_c(\varepsilon^m)$ is a positive integer
congruent to $p-m$ modulo $l$, and zero otherwise.  

(ii) If $L_c(\eta^p)\ne M_c(\eta^p)$ then $L_c(\eta^p)$ is 
finite dimensional, and the character of 
$L_c(\eta^p)$ is 
$$
\Tr(s^jt^{\bold h})=\varepsilon^{pj}
t^{\frac{1}{2}-f(\varepsilon^p)}
\frac{1-t^b\varepsilon^{-bj}}{1-t\varepsilon^{-j}},
$$
where $b$ is the smallest positive integer of the form 
$f_c(\varepsilon^p)-f_c(\varepsilon^m)$
congruent to $p-m$ modulo $l$. 
\end{theorem}

\subsection{The Gorenstein property}

Any submodule $J$ of the polynomial 
representation $M_c=\Bbb C[\h]$ is an ideal in $\Bbb C[\h]$, 
so the quotient $A=M_c/J$ is a 
$\Bbb Z_+$-graded commutative algebra. 

Now suppose that $W$ is a real reflection group. 
Recall that it was shown in \cite{BEG}, Proposition 1.13, 
that if $A$ is irreducible (i.e. $A=L_c$, the irreducible quotient of $M_c$) 
then $A$ is a Gorenstein algebra (see \cite{E} pp.529, 532 for 
definitions). 
Here we prove the converse statement.

\begin{theorem} \label{Gorirr} If 
$A=M_c/J$ is finite dimensional and Gorenstein, then $A=L_c$
(i.e. $A$ is irreducible). 
\end{theorem}

\begin{proof}
Since $A$ is Gorenstein, the highest degree component of $A$ is 
1-dimensional, and the pairing $E: A\otimes A\to \Bbb C$ given by 
$E(a,b):= {\rm h.c.}(ab)$ (where ${\rm h.c.}$ stands for the 
highest degree coefficient) is nondegenerate. 
This pairing obviously satisfies the condition 
$E(xa,b)=E(a,xb), x\in \h^*$. Now recall that $H_c(W)$ 
and $A$ admits a natural action of the group $SL_2(\Bbb C)$. Let 
$g=\begin{pmatrix} 0&1\\ -1&0\end{pmatrix}\in SL_2(\Bbb C)$; 
then $g(x_i)=y_i$ for orthonormal bases $x_i$, $y_i$ of $\h^*$ and
$\h$ which are dual to each other. Thus the nondegenerate form 
$\tilde B(a,b):=E(a,gb)$ on $A$ satisfies the 
equations $\tilde B(a,x_ib)=\tilde B(y_ia,b)$.
So for any $f_1,f_2\in \Bbb C[\h]$, 
one has $\tilde B(f_1(x)v,f_2(x)v)=
\tilde B(f_2(y)f_1(x)v,v)$, where $v=1$ is the lowest weight vector
of $A$. This shows that $\widetilde B$ coincides with the Shapovalov form $B$ 
on $A$. Thus $A$ is an irreducible representation of $H_c(W)$. 
\end{proof} 

{\bf Remark.} It easy to see by considering the rank 1 case that 
for complex reflection groups Theorem \ref{Gorirr} is, in
general, false.  

\begin{theorem} \label{fingor}
Let $W$ be a complex reflection group, and $U\subset M_c$ 
be a $W$-subrepresentation of dimension $\ell=\dim(\h)$ 
sitting in degree $r$, consisting of singular vectors 
(i.e. those killed by $y\in \h$). 
Let $J$ be the ideal generated by $U$. Assume that the
quotient representation $A=M_c/J$ is finite dimensional. 
Then 

(i) The algebra $A$ is Gorenstein. 

(ii) The representation $A$ admits a BGG type resolution 
$$
A\leftarrow M_c(\Bbb C)\leftarrow M_c(U)\leftarrow M_c(\wedge^2U)
\leftarrow...\leftarrow M_c(\wedge^\ell U)\leftarrow 0. 
$$

(iii) The character of $A$ is given by the formula 
$$
\chi_A(g,t)=t^{\frac{\ell}{2}-\sum_s \frac{2c_s}{1-\lambda_s}}
\frac{\det|_U(1-gt^r)}{\det_{\h^*}(1-gt)}.
$$
In particular, the dimension of $A$ is $r^\ell$. 

(iv) If $W$ is a real reflection group, then $A$ is irreducible. 
\end{theorem}

\begin{proof}
(i) Since ${\rm Spec}(A)$ is a complete intersection,
$A$ is Gorenstein (\cite{E} p.541).

(ii) Consider the subring $\Bbb C[U]$ in $\Bbb C[\h]$.  
Then $\Bbb C[\h]$ is a finitely generated $\Bbb C[U]$-module, 
and therefore by a standard theorem of Serre \cite{S}, it is 
free, of rank is $r^\ell$ (the latter follows from computing the 
Hilbert series). Consider the Koszul resolution attached to this module. 
At the level of $\Bbb C[\h]$ modules, this resolution 
looks exactly as we want in (ii). So we need to show 
that the maps of the resolution are in fact morphisms 
of $H_c(W)$-modules and not only $\Bbb C[\h]$-modules.
This is easily established by induction (going from left to right),
cf. proof of Proposition 2.2 in \cite{BEG} and also \cite{Go}. 

(iii) Follows from (ii) by the Euler-Poincare principle. 

(iv) Follows from Theorem \ref{Gorirr}.
\end{proof}
 
\section{Representations of the rational Cherednik algebra of type $A$}

\subsection{The results} 

Let $W=S_n$. In this case the function $c$ reduces to one number $k$. 
We will denote the rational Cherednik algebra $H_k(S_n)$ by $H_k(n)$. 
The polynomial representation $M_k$ of this algebra is the space of 
$\Bbb C[x_1,...,x_n]^T$ of polynomials of $x_1,...,x_n$, which are 
invariant under simultaneous translation $x_i\mapsto x_i+a$. 
In other words, it is the space of regular functions
on $\h=\Bbb C^n/\Delta$, where $\Delta$ is the diagonal.

\begin{proposition} \cite{Du}\label{Asing} 
Let $r$ be a positive integer not divisible by $n$, and $k=r/n$. 
Then $M_k$ contains a copy of the reflection representation 
$\frak h$ of $S_n$, which consists of singular vectors. 
This copy sits in degree $r$ and is spanned 
by the functions 
$$
f_i(x_1,...,x_n)={\rm Res}_\infty [(z-x_1)...(z-x_n)]^{\frac{r}{n}}
\frac{dz}{z-x_i}.
$$
(the symbol ${\rm Res}_\infty$ denotes the residue 
at infinity). 
\end{proposition} 

{\bf Remark.} The space spanned by $f_i$ is $n-1$-dimensional, since 
$\sum_i f_i=0$ (this sum is the residue of an exact differential). 

\begin{proof}
This proposition can be proved by a straightforward computation. 
The functions $f_i$ are a special case of 
Jack polynomials.
\end{proof}
              
\bigskip

Let $I_k$ be the submodule of $M_k$ 
generated by $f_i$, 
Consider the $H_k(n)$-module $V_k=M_k/I_k$, and 
regard it as a $\Bbb C[\h]$-module.

Our result is 

\begin{theorem}\label{main} Let $d=(r,n)$ denote the greatest common 
divisor of $r$ and $n$. 
Then the (set-theoretical) 
support of $V_k$ is the union of $S_n$-translates 
of the subspaces of $\Bbb C^n/\Delta$,
defined by the equations 
$$
\begin{array}{c}
x_1=x_2 = \dots =x_{\frac{n}{d}};\\
x_{\frac{n}{d}+1}= \dots =x_{2\frac{n}{d}};\\
 \dots  \\
x_{(d-1)\frac{n}{d}+1}= \dots = x_n.
\end{array}
$$
In particular, the Gelfand-Kirillov dimension of $V_k$ 
is $d-1$. 
\end{theorem}

The theorem allows us to give a simple proof 
of the following result of \cite{BEG} (without the use of 
the KZ functor and Hecke algebras). 

\begin{corollary} (\cite{BEG}) If $d=1$ 
then the module $V_k:=M_k/I_k$ is finite dimensional, irreducible, 
admits a BGG type resolution, and its character is 
$$
\chi_{V_k}(g,t)=t^{(1-r)(n-1)/2}\frac{\det|_\h(1-gt^r)}{\det|_\h(1-gt)}.
$$
\end{corollary}

\begin{proof} 
For $d=1$ 
Theorem \ref{main} says that the support of 
$M_k/I_k$ is $\{0\}$. This implies that 
$M_k/I_k$ is finite dimensional. The rest follows from  
Theorem \ref{fingor}.
\end{proof} 

\subsection{Proof of Theorem \ref{main}}
The support of $V_k$ is the zero-set of $I_k$, i.e. 
the common zero set of $f_i$. 
Fix $x_1,...,x_n\in \Bbb C$. Then $f_i(x_1,...,x_n)=0$ for all
$i$ iff 
$\displaystyle \sum_{i=1}^n \lambda_i f_i=0$ 
for all $\lambda_i$, i.e.
$$
 \operatorname{Res}_{\infty}\left( 
\prod_{j=1}^n(z-x_j)^{\frac{r}{n}}\sum_{i=1}^n\frac{\lambda_i}{z-x_i} 
\right)dz=0.  
$$

Assume that $x_1, \dots x_n$ take distinct values $y_1, \dots, y_p$ with 
positive multiplicities $m_1, \dots, m_p$.
The previous equation implies that the point 
$(x_1,...,x_n)$ is in the zero set iff
$$
{\rm Res}_\infty 
\prod_{j=1}^p (z-y_j)^{m_j\frac{r}{n}-1} \left(\sum_{i=1}^p
\nu_i(z-y_1)\dots \widehat{(z-y_i)} \dots (z-y_p)\right)dz=0 \quad
\forall \nu_i.
$$
Since $\nu_i$ are arbitrary, this is equivalent to the 
condition
$$
{\rm Res}_\infty 
\prod_{j=1}^p (z-y_j)^{m_j\frac{r}{n}-1}z^i dz=0, \quad i=0, \dots, 
p-1. 
$$

We will now need the following lemma. 

\begin{lemma}\label{mainlemma}
Let $\displaystyle a(z)=\prod_{j=1}^p(z-y_j)^{\mu_j}$, where 
$\mu_j\in \Bbb C$, $\sum_j \mu_j 
\in \mathbb Z$ 
and $\sum_j \mu_j > -p$. Suppose 
$${\rm Res}_\infty a(z)z^idz=0, \quad i=0,1,\dots, p-2.$$
Then $a(z)$ is polynomial.
\end{lemma}

\medskip

\begin{proof}
Let $g$ be a polynomial.
Then 
$$0={\rm Res}_\infty d(g(z)\cdot a(z))=
{\rm Res}_\infty(g^{\prime}(z)a(z)+a^{\prime}(z)g(z))dz$$

and hence

$${\rm Res}_\infty \left(g^{\prime}(z)+ \sum_i 
\frac{\mu_j}{z-y_j}g(z)\right)a(z)dz=0.$$

\smallskip

Let $\displaystyle g(z)=z^l \prod_j(z-y_j)$. Then $ 
\displaystyle g^{\prime}(z)+ \sum 
\frac{\mu_j}{z-y_j}g(z)$ is a polynomial of degree $l+p-1$ with highest 
coefficient $l+p+\sum \mu_j \ne 0$ (as $\sum \mu_j>-p$). 
This means that for every $l \ge 0$,
$\displaystyle {\rm Res}_\infty z^{l+p-1}a(z)dz$ is linear combination of 
residues of $z^qa(z)dz$ 
with $q<l+p-1$. By the assumption of the lemma, 
this implies by induction in $l$ 
that all such residues are $0$ and hence $a$ is a polynomial.
\end{proof}

In our case $\sum (m_j \frac{r}{n}-1)=r-p$ (since $\sum m_j=n$) and 
the conditions of the lemma are satisfied. Hence 
$(x_1,...,x_n)$ is in the zero set of $I_k$ iff
$\displaystyle\prod_{j=1}^p(z-y_j)^{m_j\frac{r}{n}-1}$ is a polynomial.
This is equivalent to saying that all  $m_j$ are divisible by 
$\frac{n}{d}$. 

We have proved that $(x_1, \dots x_n)$ is in the zero set of $I_k$ iff
$(z-x_1) \dots (z-x_n)$ is the $n/d$-th power of a polynomial of degree $d$.
This implies the theorem. 

\section{Representations of the rational Cherednik algebra 
for the complex reflection group 
$S_n \ltimes \left( \BZ /l \BZ\right)^n$.}

\subsection{The formula for the singular vector.}
Let $l\ge 2$ be an integer. 
Consider the complex reflection group 
$W=S_n \ltimes \left( \BZ /l \BZ\right)^n$ acting on the $n$-dimensional 
space. The set $S \subset W$  
of complex reflections consists of the elements $s_i^m$ and 
$\sigma_{i,j}^{(m)}$ defined by
$$s_i^{m}(x_1, \dots, x_i, \dots, x_n) = 
(x_1, \dots, \eps^{m}x_i,\dots, x_n), \quad 1\le m < l;$$ 
$$\sigma_{i,j}^{(m)}(x_1 ,\dots, x_i ,\dots, x_j, \dots, x_n)=
(x_1, \dots, \eps^{m}x_j, \dots, \eps^{-m}x_i ,\dots, x_n) 
,\quad 0 \le m <l
$$
(here $\eps$ is a primitive $l$-th root of unity).
Consider the $W$-invariant function $c:S \to \BC$ defined by 
$c(s_i^{m})=c_{m}$ for all $i$ and $c(\sigma_{i,j}^{(m)})=k$ for all $i, j, 
m$ (where $k$ and $\{c_m\}$ are fixed constants). Let
$H_{c}(n,l)$ 
denote the rational 
Cherednik algebra $H_c(W)$ 
corresponding to $W$ and $c$. Let $M_c$ be the polynomial 
representation of this algebra. 

Fix a positive integer $r$, which is not divisible by $l$. 
Thus $r=(p-1)l+q$, where $p$ is a positive integer, and 
$1\le q\le l-1$ is an integer. 
Denote by $E_r$ the affine hyperplane of those functions $c$ 
for which 
$l(n-1)k+2\sum_{j=1}^{l-1}c_j\frac{1-\eps^{-jq}}{1-\varepsilon^{-j}}=r$.
Let $\h_q$ be the representation of 
$W$ on $\Bbb C^n$ in which $S_n$ acts by permutations, and 
$s_i$ multiplies the $i$-th coordinate of a vector by
$\varepsilon^{-q}$
(thus the reflection representation corresponds to $q=l-1$ and
the dual reflection representation to $q=1$). 
Let $s$ be the largest integer which is $< p/n$.

We have the following analog of Proposition \ref{Asing}.

\begin{proposition}
For $c\in E_r$ the polynomial 
representation $M_c$ of $H_c(n,l)$
contains a copy of representation $\h_q$ of $W$ in degree
$r$ consisting of 
singular vectors. This copy is spanned by 
the functions $f_i$, where 
$$f_i(x_1, \dots, x_n)= 
\operatorname{Res}_{\infty}\frac{z^{(p-nk)l-1}}{(k-1)\dots(k-s)} 
\left(\prod_{j=1}^n(z^l-x_j^l)^k\right)\frac{x_i^q dz}{z^l-x_i^l} .$$
\end{proposition}

{\bf Remark.} Clearly, if $k=1,..,s$ then 
$${\rm Res}_\infty z^{(p-nk)l-1}
\left(\prod_{j=1}^n(z^l-x_j^l)^k\right)\frac{x_i^qdz}{z^l-x_i^l}=0.$$ Thus 
the functions $f_i$ are polynomial in $k$. Moreover, 
they do not vanish identically for any $k$, 
as easily seen by computing the coefficients.

\begin{proof} 
This proposition can be proved by direct computation. 
It can also be obtained as a 
simple consequence of the results of \cite{DO}.
\end{proof} 

\subsection{Finite dimensional representations of 
$H_{c}(n,l)$.}

Let $U_r=\Bbb C[u]/(u^r)$. This is naturally a
$\Bbb Z_+$-graded representation of $\Bbb Z/l\Bbb Z$, 
$m\circ u^j=\varepsilon^{-mj}u^j$. 
Thus $U_r^{\otimes n}$ is a graded representation of $W$. 
The character of this representation is 
${\rm Tr}|_{U_r}(gt^D)=
\frac{\det|_{\h_q}(1-t^rg)}{\det|_{\h^*}(1-tg)}$,
where $D$ is the grading operator. 

\begin{theorem}
(i) For $c\in E_r$, there exists a lowest weight 
module $Y_c$ over $H_{c}(n,l)$
with trivial lowest weight (i.e. a quotient of
$M_c$) which is isomorphic to $U_r^{\otimes n}$ as a graded
$W$-module.

(ii) For generic $c\in E_r$, $Y_c$ is irreducible.  
\end{theorem}
 
\begin{proof} Let $I_c$ be the ideal generated by 
$f_i$, and $\widetilde{Y}_c=M_c/I_c$. 

If $k=0$, a direct computation shows that 
$f_i(x_1,...,x_n)=Cx_i^r$, where $C$ is a nonzero constant.
In this case, $H_{c}(n,l)=\Bbb C[S_n]\ltimes 
(H_{c_1,...,c_{l-1}}(1,l)^{\otimes n})$, so
$\widetilde{Y}_c$ is simply $U_r^{\otimes n}$, where $U_r$ is 
the $r$-dimensional lowest weight module 
over $H_{c_1,...,c_{l-1}}(1,l)$ (with trivial lowest weight) which
exists when $2\sum_{j=1}^{l-1}c_j\frac{1-\eps^{-jq}}{1-\varepsilon^{-j}}=r$.  
It follows from Theorem \ref{rank1} that 
for generic $c_1,...,c_{l-1}$ satisfying this equation, 
the module $U_r$ is irreducible. Thus, for generic $c\in E_r\cap
\lbrace{k=0\rbrace}$, 
the module $\widetilde{Y}_c$ is irreducible (of dimension
$r^n$), and the kernel of the Shapovalov form on $M_c$
coincides with $I_c$. 

Using standard semicontinuity arguments, we conclude from this that for
generic $c\in E_r$, ${\rm dim}(\widetilde{Y}_c)\le r^n$. 

On the other hand, let $L_c$ denote the irreducible module 
over $H_c(n,l)$ with trivial lowest
weight, i.e. the quotient of $M_c$ by the kernel 
of the Shapovalov form. The previous argument shows 
that for generic $c\in E_r$, ${\rm dim}L_c\ge r^n$. 
But $L_c$ is clearly a quotient of $\widetilde{Y}_c$. This implies that for
generic $c\in E_r$, $L_c=\widetilde{Y}_c$, and in particular 
$\widetilde{Y}_c$ is
irreducible. 

To finish the proof, observe that the ideal $I_c$
really depends only on one parameter $k$. 
Therefore, we can define $Y_c$ by setting, 
for generic $k$, $Y_c:=\widetilde{Y_c}=M_c/I_c$, and for special $k$
(i.e. finitely many values), $Y_c=M_c/\hat I_c$, where $\hat
I_c=\lim_{c'\to c}I_{c'}$ (the limit exists since 
any regular map from $\Bbb P^1\setminus \lbrace{k_1,...,k_N\rbrace}$ 
to a Grassmannian can 
be extended to the whole $\Bbb P^1$). Clearly, 
$Y_c$ is a lowest weight $H_c(n,l)$-module, satisfying 
the conditions of the theorem. Thus the theorem is proved.    
\end{proof}

{\bf Remark.} Let $l=2$. In this case our condition on $r$ is
that $r$ is odd, and the equation of $E_r$ has the form
$2(n-1)k+2c_1=r$. Thus, in this case we recover Theorem 6.1 from
\cite{BEG}.

\subsection{The dimension of $\widetilde{Y_c}$}

Let 
$$
\Sigma_r=\lbrace{\frac{P}{Q}| \quad P, Q \text{ 
integers and}\quad (P,Q)=1,\quad
1\le P\le p-1, \quad 1 \le Q\le n \rbrace}
$$ 
(where 
$r=(p-1)l+q$), and let $\widetilde{Y_c},\ Y_c$ be as defined in
 the previous subsection. 

\begin{theorem}\label{fini} (i) $\widetilde{Y_c}$ is finite
dimensional if and only if $k\notin \Sigma_r$. 

(ii) $\widetilde{Y_c}$ is finite
dimensional if and only if $\widetilde Y_c=Y_c$. 

(iii) If $l=2$ and $k\notin \Sigma_r$ then 
$Y_c$ is irreducible. 
\end{theorem}

\begin{proof}
(i) It suffices to prove that the support of $\widetilde{Y_c}$ is $\{0\}$ iff 
$k \notin \Sigma_r$. The 
proof will be analogous to that of Theorem \ref{main}.

The support of $\widetilde{Y_c}$ is the common zero set of the
functions $f_i$. 
Fix $x_1,...,x_n\in \Bbb C$. 
We have: $f_i(x_1, \dots, x_n)=0$ for all $i$ iff $\displaystyle 
\sum_{i=1}^n \lambda_i f_i(x_1,...,x_n) =0$ for all $\lambda_i$, 
i.e.
$$
{\rm Res}_{\infty} \left( z^a \prod_{j=1}^n(z^l-x_j^l)^k 
\sum_{i=1}^n \frac{\lambda_i x_i^q}{z^l-x_i^l}\right)dz=0,
$$
where $a=(p-nk)l$. 

Let $d\ge 0$ be the number of distinct nonzero  
numbers among $x_1^l,...,x_n^l$. More specifically, 
assume that $x_1^l, \dots, x_n^l$ take  
values $y_0=0,y_1, \dots, y_d$ with multiplicities $m_0,m_1,
\dots, m_d$,
such that $m_j> 0$ if $j>0$ (so $\sum_{i=0}^d m_j=n$).
From the above equation we see that 
the point $(x_1,...,x_n)$ is in the zero set if and only if
for any $\nu_i$, $i=1,...,d$
$$
{\rm Res}_\infty\left[ z^{a+lm_0} 
\prod_{j=1}^d(z^l-y_j)^{m_jk-1}\left(\sum_{i=1}^{d}
\nu_i (z^l-y_1) \dots \widehat{(z^l-y_i)} \dots (z^l-y_d) 
\right) \right]z^{-1}dz=0.
$$ 
This equation is equivalent to
$$
{\rm Res}_\infty \left(z^{a+(m_0k+i)l} 
\prod_{j=1}^d(z^l-y_j)^{m_jk-1}\right)z^{-1}dz =0, 
\quad i=0, \dots, d-1.
$$
Lets make a change of variables $w=z^l$.
Then the equations above will take the form
$$
{\rm Res}_\infty 
\left((w-y_0)^{p-(n-m_0)k-1}
\prod_{j=1}^d(w-y_j)^{m_jk-1}w^i\right)dw=0, \quad 
i=0,\dots,d-1.
$$ 
Applying Lemma \ref{mainlemma} 
(whose conditions are clearly satisfied) we get that
these equations hold iff the function
$$
F(w):=w^{p-(n-m_0)k-1}
\prod_{j=1}^d(w-y_j)^{m_jk-1}
$$
is a polynomial.

The function $F(w)$ is a polynomial iff  $p-(n-m_0)k-1$ and  $m_j k-1$ for all 
$1 \le j \le d$  
are non-negative integers.

Suppose $k \notin \mathbb Q $ or $k=\frac{P}{Q}$ with 
$(P,Q)=1$ and $Q > n$. Then there 
is no integer $1 \le m \le n$, such that $mk-1 \in \mathbb Z$. This means 
that for such $k$ 
$F$ is a polynomial only 
if $d=0$. Hence the support of $\widetilde{Y_c}$ is $\{0\}$.

Suppose $k=\frac{P}{Q}$ with $(P,Q)=1$ and $Q \le n$.
Let $M_j=km_j$. The condition that $F$ is a polynomial 
implies that $M_j$ are positive integers and $p-1-\sum_{j=1}^d
M_j\ge 0$. This means that either $d=0$ or $P\le p-1$, i.e. $k\in
\Sigma_r$. Thus, if the support of $\widetilde{Y_c}$
is nonzero, then $k\in \Sigma_r$. 
 
Conversely, let $k=P/Q\in \Sigma_r$. Take $d=1$, $m_1=Q$, 
$m_0=n-Q$, and choose $y_1\ne 0$ arbitrarily. Then $F$ is a polynomial. So 
the support of $\widetilde{Y}_c$ in this case is nonzero. 

Thus, statement (i) is proved.  

(ii) Clearly, $Y_c$ is a quotient of $\widetilde{Y_c}$. 
The rest follows from Theorem \ref{fingor}. 

(iii) If $l=2$ then $W$ is a reflection group. 
Thus the result follows from Theorem \ref{fingor}.
\end{proof} 

Part (iii) of Theorem \ref{fini} generalizes Proposition 6.4 of
\cite{BEG}. 

{\bf Remark.} Part (iii) of Theorem \ref{fini} fails for 
complex reflection groups, as seen from considering the rank 1
case: in this case $\Sigma_r$ is empty, but $Y_c$ is not always
irreducible. Similarly, one cannot drop in part (iii) the
assumption that $k\notin \Sigma_r$: this is demonstrated by
Subsection 6.4 in \cite{BEG}.

\end{document}